\documentclass[11pt]{article}

\usepackage{authblk}
\usepackage{amsfonts}
\usepackage{amsthm}
\usepackage{amsmath}
\usepackage{graphicx}
\usepackage{comment}
\usepackage[usenames,dvipsnames]{xcolor}
\usepackage{enumerate}
\usepackage{extarrows}
\usepackage{hyperref}
\usepackage[round, authoryear]{natbib}
\usepackage{stackrel}
\usepackage{centernot}
\usepackage{caption}
\usepackage{subcaption}
\usepackage{fancyhdr}

\let\quoteOLD\quote
\def\quote{\quoteOLD\small}

\definecolor{labelkey}{cmyk}{0,0.8,1,0.5}
\definecolor{refkey}{cmyk}{0,0.8,1,0.5}

\newtheorem{theorem}{Theorem}
\newtheorem{corollary}{Corollary}
\newtheorem{proposition}{Proposition}
\newtheorem{definition}{Definition}
\newtheorem{lemma}{Lemma}
\newtheorem{remark}{Remark}

\numberwithin{equation}{section}
\numberwithin{theorem}{section}
\numberwithin{corollary}{section}
\numberwithin{proposition}{section}
\numberwithin{lemma}{section}
\numberwithin{definition}{section}
\numberwithin{remark}{section}

\newtheorem{example0}{\sc Example}[subsection]

\makeatletter
\def\th@newremark{\th@remark\thm@headfont{\bfseries}}
\makeatletter

\def\nexto{\kern -0.54em}

\def\boxit#1{\vbox{\hrule\hbox{\vrule\kern6pt 
          \vbox{\kern6pt#1\kern6pt}\kern6pt\vrule}\hrule}}

\newcommand{\laminv}{\overline{\Lambda}^{\leftarrow}}

\newcommand{\eqd}{\stackrel{{\mathrm D}}=}

\newcommand{\AAA}{{\cal A}}\newcommand{\BBB}{{\cal B}}

\newcommand{\FFF}{{\cal F}}

\newcommand{\JJJ}{{\cal J}}
\newcommand{\MMM}{{\cal M}}

\newcommand{\BB}{\mathbb{B}}

\newcommand{\MM}{\mathbb{M}}\newcommand{\NN}{\mathbb{N}}

\newcommand{\XX}{\mathbb{X}}

\newcommand{\R}{\Bbb{R}}

\newcommand{\N}{\Bbb{N}}

\newcommand{\wt}{\tilde}

 \newcommand{\EEEE}{\mathfrak{E}}

\newcommand{\be}{\begin{equation}}\newcommand{\ee}{\end{equation}}
\newcommand{\bea}{\begin{eqnarray}}\newcommand{\eea}{\end{eqnarray}}
\newcommand{\bean}{\begin{eqnarray*}}\newcommand{\eean}{\end{eqnarray*}}
\newcommand{\ben}{\begin{equation*}}\newcommand{\een}{\end{equation*}}
\newcommand{\ba}{\begin{aligned}}\newcommand{\ea}{\end{aligned}}

\newcommand{\halmos}{\quad\hfill\mbox{$\Box$}}

\newcommand{\BN}{\mathcal{BN}}
\newcommand{\NB}{\mathrm{NB}}
\newcommand{\PPP}{\mathrm{PPP}}

\newcommand{\bfV}{{\bf V}}

\newcommand{\bfW}{{\bf W}}

\newcommand{\PD}{\textbf{\rm PD}}
\newcommand{\PK}{\textbf{\rm PK}}

\newcommand{\rmd}{{\rm d}}

\newcommand{\PP}{\mathbb{P}}
\newcommand{\EE}{\mathbb{E}}

\begin{document}

\title{ Negative Binomial Construction of Random Discrete Distributions on the Infinite Simplex}

\author{Yuguang F. Ipsen\thanks{Corresponding author: Yuguang.Ipsen@anu.edu.au}  }
\author{Ross A. Maller\thanks{Research partially supported by ARC Grant DP1092502; Email: Ross.Maller@anu.edu.au}}
\affil{Research School of Finance,  Actuarial Studies \&  Statistics, \\
Australian National University, Canberra, Australia.}

\maketitle

\begin{abstract}
The Poisson-Kingman distributions,  $\PK(\rho)$,  on the infinite simplex, can be constructed from a Poisson point process having intensity density $\rho$
or by taking the ranked jumps up till a specified time   
of a subordinator
 with L\'evy density $\rho$, 
 as proportions of the subordinator. 
 As a natural extension, we replace the Poisson point process with a negative binomial point process having parameter $r>0$ and 
 L\'evy density $\rho$,  thereby defining a new class
 $\PK^{(r)}(\rho)$ 
 of distributions on the infinite simplex.
The new class contains the  two-parameter generalisation $\PD(\alpha, \theta)$ of \cite{PY1997} when $\theta>0$.
It also contains a class of distributions derived from the trimmed stable subordinator.
We derive properties of the new distributions,
with particular reference to the two most well-known 
 PK distributions: the Poisson-Dirichlet  distribution $\PK(\rho_\theta)$ 
 generated by a Gamma process with L\'evy density  $\rho_\theta(x) = \theta e^{-x}/x$, $x>0$, $\theta > 0$, and the random discrete distribution,  $\PD(\alpha,0)$, 
derived from an $\alpha$-stable subordinator.

\vspace{0.3cm}


\end{abstract}

\section{Introduction}
The random discrete distributions on the infinite simplex constructed from Poisson point processes possess both elegant theoretical properties and wide applicability to many real-world problems. The construction is as follows. Let $\XX = \sum_{i=1}^\infty \delta_{\Delta_{(i)}}$,  where $\delta_x$ denotes a point mass at $x\in\R$,  be a Poisson point process on $(\R_+, \BBB(\R_+))$ with ordered points
 $\Delta_{(1)} \ge \Delta_{(2)} \ge \cdots$ and intensity measure $\Pi(\cdot)$, satisfying
\be\label{sub}
\Pi\{(0,\infty)\} = \infty, \,\, \Pi\{(x,\infty)\} < \infty \,\, \text{for each } x > 0, \,\, \text{and} \,\, \int_{0}^1 x \, \Pi(\rmd x) < \infty. 
\ee
Denote the sum of points in $\XX$ by $T(\XX) = \sum_{i} \Delta_{(i)} $. Then condition \eqref{sub} ensures that $\PP(0 < T(\XX) < \infty) = 1$.  Assume 
that $\Pi$ admits a density $\rho$, so that $\Pi(\rmd x) = \rho(x)\rmd x$ for any $x >0$.
Define $\bfV$ to be the ordered jumps normalised by their sum:
\be\label{PD}
\bfV = (V_1, V_2, \ldots) = \bigg(\frac{\Delta_{(1)}}{T(\XX)}, \,   \frac{\Delta_{(2)}}{T(\XX)}, \ldots \,  \bigg).
\ee
Then $\bfV$ follows a Poisson-Kingman distribution with density $\rho$, denoted by $\PK(\rho)$, following the terminology and notation in \cite{Pitman2003}. 
The two most well-known Poisson-Kingman distributions are the Poisson-Dirichlet distribution, $\PK(\rho_\theta)$,  where $\rho_\theta(x) = \theta e^{-x}/x$, $x>0$, $\theta > 0$, introduced in \cite{kingman1975}, and the random discrete distribution derived from an $\alpha$-stable subordinator, $\PK(\rho_\alpha)$, with $\rho_\alpha(x) = C \alpha x^{-\alpha-1}$,   $x>0$, for some $C >0$ and $0<\alpha<1$. \footnote{We will always use the Greek letters $\theta$ and $\alpha$, and only these, to parametrise the $\rho_\theta$ and $\rho_\alpha$ as just defined, so there should be no confusion between them.}
This distribution is also known as $\PD(\alpha,0)$ in \cite{PY1997}, in which a two-parameter generalisation $\PD(\alpha, \theta)$ is constructed. These distributions have applications ranging from the modelling of gene frequencies in population genetics through to models for prior distributions in Bayesian nonparametric statistics as well as, recently in the machine learning community;  for example,  \cite{Ewens1972};  \cite{Donnelly1986}; \cite{James2001};  \cite{Lim2016}.

The construction of $\PK(\rho)$ from a Poisson point process is equivalent to that from a subordinator. Let $X=(X_t)_{t\ge 0}$ be a driftless subordinator with L\'evy density $\rho$. 
Write the jump process of $X$ as $(\Delta X_t := X_t- X_{t-})_{t>0}$, and its jumps up till time $t>0$ arranged in decreasing order as $\Delta X_t^{(1)} \ge \Delta X_t^{(2)}\ge \cdots $. 
Then by the definition in \eqref{PD}, the vector
\be\label{csub}
\bigg(\frac{\Delta X_t^{(1)}}{X_t}, \,  \frac{\Delta X_t^{(2)}}{X_t}, \cdots \,  \bigg) \sim \PK(t\rho),
\ee
where ``$\sim$" stands for ``has the distribution of".

As a natural extension of \eqref{csub}, one can consider the subordinator $X$ after omitting a fixed number, $r\in\N$,  of its largest jumps, and forming the analogous normalised vector of ratios.
The corresponding {\it $r$-trimmed 
subordinator} is ${}^{(r)}X_t := X_t-
\sum_{i=1}^r \Delta X_t^{(i)}$.
When $X$ is an $\alpha$-stable subordinator, we denote the generalised class of distributions by $\PD_\alpha^{(r)}$, so that
\be\label{PDr}
\bigg(\frac{\Delta X_1^{(r+1)}}{{}^{(r)}X_1}, \,  \frac{\Delta X_1^{(r+2)}}{{}^{(r)}X_1}, \cdots \,  \bigg) \sim \PD_\alpha^{(r)}\,.
\ee
There is a close connection between the trimmed stable subordinator and the negative binomial point process of \cite{Gregoire1984} which we develop in Section~\ref{sub:trimmed}.
The ``trimming" concept, of removing the $r$ largest points, is very natural in this context, but having defined the new distributions, there is no need to keep $r$ as an integer, and any value $r>0$ can be allowed (as we do herein).

Returning to the general situation, this suggests replacing
the Poisson construction with a negative binomial construction.
Alternatively to the trimming rationale mentioned in the previous paragraph, this generalisation is in the same spirit as replacing a Poisson by a negative binomial distribution in the statistical analysis of discrete data, to cater for overdispersion in the data.
In any event, this procedure produces a new class of random discrete distributions, parametrised by an extra parameter $r>0$, which includes the important special case $\PD(\alpha, \theta)$ for $\theta > 0$, hence, in particular, the classes $\PD(\alpha,0)$ and $\PD(0,\theta)$, as well as $\PD_\alpha^{(r)}$. We go on to study stick-breaking representations and other properties of the new class.

The paper is organised as follows. Section \ref{sect:NB}  sets up the negative binomial construction and defines the new Poisson-Kingman class, denoted by $\PK^{(r)}(\rho)$. We note the connection of this class to subordinated L\'evy processes as well as to the $\PD(\alpha, \theta)$ distribution when $\theta > 0$.
From previous Poisson formulae, we then derive the joint density of the size-biased random permutation in Section \ref{sub:ssb} and note its interesting consequences. 
In Section \ref{sub:stick}, we give the stick-breaking representation for important cases of the new class (see Theorem \ref{stick}) and make explicit how it differs from the original $\alpha$-stable case. 
Section \ref{sub:trimmed} defines a special case $\PD_\alpha^{(r)}$ of the new family that arises naturally from trimmed stable subordinators as discussed above.
The related stick-breaking property is investigated in Theorem \ref{stick2} and a characterisation as a shifted sequence from $\PD(\alpha, 0)$ is obtained in Lemma \ref{changeM}. 
Finally, the paper ends with a brief discussion of potential applications of the new class in Section \ref{dis}.

\section{The Negative Binomial Construction}\label{sect:NB}
Throughout, we assume a background probability space $(\Omega, \AAA, \PP)$ is given and all random variables are measurable mappings from $(\Omega, \AAA, \PP)$ to some appropriate space. 
We follow the exposition in \citet[Ch.3]{res87} for the point process setup. Denote the space of integer-valued Radon measures on $(\R_+, \BBB(\R_+))$ endowed with its Borel $\sigma$-algebra by $(\MM, \MMM)$, and let $\FFF^+$ be the set of nonnegative measurable functions on $(0,\infty)$. 
We use the abbreviation $\PPP(\Pi)$
(or $\PPP(\rho)$) throughout to denote a 
 Poisson point process with intensity measure 
$\Pi$ (or density $\rho$).

Given a measure $\Pi$ on $\BBB(\R_+)$, locally finite at infinity,  and  a parameter $r > 0$,
a negative binomial point process (NBPP) on $(\R_+, \BBB(\R_+))$
generated by  $\Pi$ and $r$, as defined in \cite{Gregoire1984}, is a measurable mapping $\BB^{(r)}$ from $(\Omega, \AAA, \PP)$ to  $(\MM, \MMM)$ characterised by its Laplace functional
\be\label{NBp}
 \Phi(f) = \EE\Big(e^{-\BB^{(r)}(f)} \Big) =\Big(1 + \int_{\R_+}(1-e^{-f(x)})\Pi(\rmd x) \Big)^{-r}, \ f \in \FFF^+.
\ee
Denote the distribution of such a $\BB^{(r)}$  by $\BN(r, \Pi)$ or $\BN(r, \rho)$ if $\Pi$ admits a density $\rho$. 
From \eqref{NBp},  $\BB^{(r)}$ with distribution $\BN(r, \Pi)$ can be regarded as a Poisson point process with randomised intensity measure $\Gamma_r \Pi$, where $\Gamma_r$ is an independent Gamma($r$,1) random variable.

The ``negative binomial" terminology arises as follows. 
 Let $\BB^{(r)}$ be a point process on $\R_+$ distributed as $\BN(r, \Pi)$ for some $r \in \NN$,  and let $B_1, \ldots, B_n$ be a sequence of pairwise disjoint bounded Borel sets on $\R_+$. Then the numbers of points of  $\BN(r, \Pi)$ 
 in $B_1, \ldots, B_n$ follow an $n$-variate negative binomial distribution $\NB(r, q_0, q_1, \ldots, q_n)$, with
\[ q_0 =  \frac{1}{1 + \Pi\big\{\bigcup_{i=1}^n B_i\big\} } \quad \text{and} \quad 
q_i =  \frac{\Pi(B_i)}{1 + \Pi\big\{\bigcup_{i=1}^n B_i\big\}} \,\, \text{ for } i=1, \ldots, n.
\]
Write $\BB^{(r)}$ in the form
\ben
\BB^{(r)} = \sum_{i\in\N} \delta_{J_{(i)}},
\een
where $J_{(1)} \ge J_{(2)} \ge \cdots$ are the ordered points in  $\BB^{(r)}$, and denote by $T(\BB^{(r)}) := \sum_{i\in\N} J_{(i)}$ the sum of the points in $\BB^{(r)}$.
Since $\Pi$ satisfies \eqref{sub}, $\PP( 0< T(\BB^{(r)}) < \infty ) = 1$ for all $r > 0$.
We now define  from $\BN(r, \Pi)$ a random discrete distribution on the infinite simplex using a similar procedure as in \eqref{PD}.

\begin{definition}[$\PK^{(r)}(\rho)$]\label{NBD}
For each $r > 0$ let $\BB^{(r)}$ be distributed as $\BN(r, \Pi)$ and suppose $\Pi$ admits a density $\rho$.  The vector
\be\label{MPD}
\bfW := (W_1, W_2, \ldots, ) = \bigg(\frac{J_{(1)}}{T(\BB^{(r)})}, \frac{J_{(2)}}{T(\BB^{(r)})}, \ldots \bigg), 
\ee
is said to follow a Poisson-Kingman distribution generated by 
$\BN(r, \rho)$, which we  denote as $\PK^{(r)}(\rho)$.
\end{definition}

Just as $\BN(r, \rho)$ can be characterised as a $\PPP(\Gamma_r \rho)$, so 
an equivalent construction of the $\PK^{(r)}(\rho)$ sequence can be made from Gamma subordinated L\'evy processes.
Let $X=(X_t)_{t > 0}$ be a driftless subordinator with L\'evy density $\rho$ and let $(\sigma_r)_{r > 0}$ be an independent gamma subordinator, i.e., a subordinator having 
L\'evy density $\Pi_{\sigma}(\rmd z)=e^{-z}z^{-1}\rmd z$, $z>0$. Denote the ranked jumps of $X$ up till time $t>0$ by $\Delta X_t^{(1)}> \Delta X_t^{(2)}>\cdots $. 
From \eqref{csub}  we immediately obtain
\be\label{subPK}
\bigg(\frac{\Delta X_{\sigma_r}^{(1)}}{X_{\sigma_r}}, \frac{\Delta X_{\sigma_r}^{(2)}}{X_{\sigma_r}}, \frac{\Delta X_{\sigma_r}^{(3)}}{X_{\sigma_r}}, \ldots \bigg)\sim \PK^{(r)}(\rho).
\ee  

One may generalise this further by considering $\PPP(\xi \rho)$ with any positive random variate $\xi$ replacing $\Gamma_r$. 
See for example \cite{JOT} for a generalisation of the Indian Buffet process using scaled subordinators.
Our emphasis on the negative binomial class is due to its natural derivation from an $r$-trimmed stable subordinator (see Section  \ref{sub:trimmed}), but it also turns out that
 $\PK^{(r)}(\rho)$  defined in \eqref{MPD} includes the well-known two parameter Poisson-Dirichlet distributions $\PD(\alpha,\theta)$, when $\theta > 0$.

 Class $\PD(\alpha, \theta)$ is defined in \cite{PY1997} through a stick-breaking representation of size-biased permutations of vectors having those distributions. 
Thus, 
$\bfV = (V_i)_{i\ge1}$ is said to have law $\PD(\alpha, \theta)$ with size-biased permutation $\wt \bfV = (\wt V_i)_{i\ge1}$ if for each $n \in \N$
\ben 
\wt V_n = U_n \prod_{i=1}^{n-1} (1- U_i),
\een
where the $U_i$ are independent and distributed as Beta($\theta + i\alpha, 1-\alpha)$.
When $\theta >0$,  $\PD(\alpha, \theta)$ has a subordinator representation in terms of the generalised Gamma subordinator. Take $0<\alpha<1$ and denote the L\'evy density of the generalised Gamma subordinator as
\ben\label{gg}
\rho_G (x) = \frac{\alpha}{\Gamma(1-\alpha)} x^{-\alpha-1} e^{-x}, \quad x > 0\,.
\een
Let $(X_t)_{t\ge 0}$ be a subordinator with L\'evy density $\rho_G$ and let
$(\sigma_r)_{r > 0}$ be an independent gamma subordinator. 
Then, by \citet[Prop.21]{PY1997}, 
\be\label{21PY}
\bfV = \Big( \frac{\Delta X_{T}^{(1)}}{X_T},  \frac{\Delta X_{T}^{(2)}}{X_T}, \ldots    \Big)  \sim \PD(\alpha, \theta) \quad \text{if}\quad T =\sigma_{\theta/\alpha}.
\ee
Comparing \eqref{subPK} with \eqref{21PY}, we see that 
for each $\theta > 0$, $\PK^{(\theta/\alpha)}(\rho_G)$ has the same law as $\PD(\alpha, \theta)$. In particular, $X_T$ is independent of $\bfV$ in \eqref{21PY}.

\begin{remark}{\rm 
It is possible to enlarge the current class $\PK^{(r)}(\cdot)$ to include $\PD(\alpha, \theta)$ when $0 > \theta > -\alpha$ by using a  2-variate mixing measure $\eta_2(\cdot, \cdot)$ on $(0,\infty) \times (0, \infty)$ such that
\ben\label{MPK}
\PK(\rho, \eta_2(\cdot, \cdot)) = \int_{s > 0} \int_{t > 0} \PK(s\rho \,| \, t) \, \eta_2(\rmd s, \rmd t).
\een
Then the mixed class $\PK(\rho, \eta(\cdot))$ defined in  Definition 3 of \cite{Pitman2003} is included in the enlarged class with mixing measure $\eta_2(\rmd s, \rmd t) =  \delta_{1}(ds) \, \eta(\rmd t)$ trivially. So we have a generalisation of $\PD(\alpha, \theta)$ for the entire parameter range. }
 \end{remark}

In the following sections, many interesting properties of $\PK^{(r)}(\rho)$ are derived through the corresponding Poisson formulae in \cite{PPY1992} and \cite{PY1997}.

\section{Joint Density of the Size-biased Permutation}\label{sub:ssb}
The {\it size-biased permutation} of a sequence $\JJJ:=(J_{(i)})_{i\in\N}$, denoted by $(\wt J_i)_{i\in\N}$, is defined as follows.
Conditional on  $\JJJ$, 
$\wt J_1$ takes value  $J_{(i)}$ 
with probability $J_{(i)}/T(\BB^{(r)})$; 
for $n\ge 1$, 
conditional on $\JJJ$ and $\{ \wt J_1, \ldots, \wt J_n \}$, $\wt J_{n+1}$ takes value  $J_{(j)} \in \JJJ\setminus\{ \wt J_1, \ldots, \wt J_n \}$ with probability 
 $J_{(j)}/\big(T(\BB^{(r)})- \sum_{i=1}^n \wt J_i \big)$. 

Let the sum of the points remaining after each size-biased pick from $\BB^{(r)}$ be the sequence $({}^{(r)}T_n)$, $n\ge 0$;
thus,  ${}^{(r)}T := {}^{(r)}T_0 := T(\BB^{(r)})$, and,  for each $n\in\N$,
\be\label{sb_T}
{}^{(r)}T_n := \, {}^{(r)}T_{n-1} - \wt J_n\,.
\ee
For each $ r > 0 $, denote the density of the random variable  $T(\mathbb{B}^{(r)})$ by
\be\label{gd}
g_{r}(t) := \PP\big( T(\BB^{(r)}) \in \rmd t\big)/\rmd t = \PP\big({}^{(r)}T \in \rmd t\big) /\rmd t, \ t>0.
\ee
By \eqref{NBp}, $g_{r}$ exists and satisfies, for each $\lambda > 0$, 
\ben\label{gdL}
\int_0^\infty e^{-\lambda x} g_{r}(x) \rmd x = 
\Big(1 + \int_{0}^{\infty} (1- e^{-\lambda x}) \Pi(\rmd x)  \Big)^{-r}.
\een
Furthermore, we always have 
\be\label{gdL2}
 g_r (t) = \int_{v > 0} f^v (t) \PP(\Gamma_r \in \rmd v )  
\ee where $f^v$ is the density of the sum of Poisson points $T(\XX)$ with L\'evy density $v\rho$.

Since $\PK^{(r)}(\rho) = \PK(\Gamma_r \rho)$, we can derive the joint densities of the sequence of remaining sums $({}^{(r)}T_i, i\ge 0)$ by applying the Poisson formulae in \cite{PPY1992}. (This is also derived from first principles through Palm characterisation in \cite{trimfrenzy}.)
Write the ascending factorial of base $r$ and order $n \in \N$ as $r^{[n]} = r(r+1)(r+2)\cdots(r+n-1)$ with  $r^{[0]} =1$.

\begin{theorem}\label{th2.1} 
Let $r > 0$ and $\Theta(x) := x \rho(x)$, $x>0$.
\begin{enumerate}[{\rm (i)}]
\item The joint density of $\left({}^{(r)}T, \, {}^{(r)}T_1, \, {}^{(r)}T_2, \ldots,  {}^{(r)}T_n\right)$ with respect to Lebesgue measure is, for $t_0 >t_1> \cdots > t_{n} > 0$  and $n \in \NN$,
\be\label{jointT0}
f(t_0, t_1, \ldots, t_n) = r^{[n]}  g_{r+n}(t_n) \prod_{i=0}^{n-1} \frac{\Theta(t_i - t_{i+1})}{t_i}.
\ee

\item $\{{}^{(r)}T, {}^{(r)}T_1, {}^{(r)}T_2,\ldots \}$ is a non-homogeneous Markov chain with transition density 
\be\label{MCr}
\PP\big( {}^{(r)}T_{n+1} \in t_1 \,\big|\,  {}^{(r)}T_{n} = t \big) = (r+n) \frac{\Theta(t-t_1)}{t} \frac{g_{r+n+1}(t_1)}{g_{r+n}(t)}\, \rmd t_1,
\ee
for each $n = 0,1,2,\ldots $, $t > 0$ and $0 <  t_1 <t$.
\end{enumerate}
\end{theorem}

\medskip\noindent {\bf Proof of Theorem \ref{th2.1}:}\
Recall in \eqref{gdL2} that $f^v$ is the density of the sum of Poisson points with L\'evy density $v\rho(x)$. By \citet[Theorem 2.1]{PPY1992}, the joint density of the corresponding remaining sum $\left(T_i,  i =0,\ldots, n\right)$ is 
\ben
  f^v(t_n) \,\,  \prod_{i=1}^n \frac{ v^n\Theta(t_{i-1} -t_i)}{t_i}.
\een
Randomising $v$ by a Gamma$(r,1)$ distribution, we see that $\left({}^{(r)}T_i,  i =0,\ldots, n\right)$ has joint density 
\begin{align*}
\int_{v > 0}  f^v(t_n) v^n  \PP(\Gamma_r\in \rmd v)  \,\,  \prod_{i=1}^n \frac{\Theta(t_{i-1} -t_i)}{t_i} .
\end{align*}
Here we note by \eqref{gdL2} that
\[\int_{v > 0}  f^v(t_n) v^n  \PP(\Gamma_r\in \rmd v) = 
\frac{\Gamma(r+n)}{\Gamma(r)}\int_{v > 0}  f^v(t_n) \PP(\Gamma_{r+n} \in \rmd v) = r^{[n]} g_{r+n}(t_n).
\]
This proves \eqref{jointT0}.
Part (ii) follows immediately from Part (i) as
\bean 
&&\, \PP\big({}^{(r)}T_{n+1} \in \rmd t_{n+1} \, \big|\, {}^{(r)}T = t_0, \, {}^{(r)}T_1= t_1, \, {}^{(r)}T_2 = t_2, \ldots,  {}^{(r)}T_n = t_n \big) \nonumber \\
&&=\,  (r+n)  \frac{\Theta(t_n-t_{n+1})}{t_n}  \frac{g_{r+n+1}(t_{n+1})}{g_{r+n}(t_n)}\rmd t_{n+1},
\eean
which does not depend on $t_0,t_1,\ldots, t_{n-1}$.
Thus \eqref{MCr} is established.
\halmos

We next write down some interesting consequences of Theorem \ref{th2.1}. 
Suppose $\bfW = (W_i)_{i\in\N}$ is defined as in \eqref{MPD} with size-biased permutation of $(W_i)$ as $(\wt W_i)$.
Let the distribution of $\bfW$ conditional on $\{T(\BB^{(r)}) = t \}$ be $\PK^{(r)}(\rho \,|\, t)$.

\begin{corollary}\label{struct}
For each  $r > 0$ and $t > 0$, the following statements hold.

{\rm (i)} The joint density of $(\wt J_1, {}^{(r)}T)$ is
\be\label{prop3.4}
\PP\big( \wt J_1 \in \rmd v, {}^{(r)}T \in \rmd t \big) = r \frac{v}{t} \, \rho(v) \, g_{r+1}(t-v) \, \rmd v \, \rmd t,\ t>0,\ 0<v< t.
\ee

{\rm (ii)}\  $g_r$ satisfies the following integral
recursion  equation:
\be\label{gr_int}
g_r(t) = r \int_0^t \rho(v)g_{r+1}(t-v) \frac{v}{t} \rmd v, \quad t > 0.
\ee

{\rm (iii)}\ The first size-biased pick from $\PK^{(r)} (\rho \, | \, t)$ has density
\be\label{con_density}
\wt f_r(w \, | \, t) =  \PP \big( \wt W_1 \in \rmd w \,\big| \, {}^{(r)}T = t \,  \big)/\rmd w =r \, tw\, \rho(tw) \, \frac{g_{r+1}(t \bar w)}{g_r(t)}, 
\ee
where $\bar w = 1 - w$ and $w \in (0,1)$.

{\rm (iv)}\ Let $G_{i+1} = \wt J_{i+1}/{}^{(r)}T_i$. For each $i = 0,1,\ldots$, $w \in (0,1)$,
\be\label{gr_MC1}
\PP \big(G_{i+1} \in \rmd w \, \big| \, {}^{(r)}T_0, \ldots,{}^{(r)}T_i = t  \, \big) 
= \PP \big(G_{i+1} \in \rmd w \, \big|{}^{(r)}T_i = t \, \big) 
= \wt f_{r+i}(\, w \,|\, t) \rmd w.
\ee
\end{corollary}

\medskip\noindent {\bf Proof of Corollary \ref{struct}.}\
(i)\ The lefthand side of \eqref{prop3.4} equals $\PP\big( {}^{(r)}T  \in \rmd t, {}^{(r)}T_1 \in \rmd (t-v) \big)$. Apply the density formulae in \eqref{jointT0} with $f(t, t-v)$ to get the righthand side of \eqref{prop3.4}.

(ii)\ Integrate \eqref{prop3.4} with respect to $\rmd v$ to obtain \eqref{gr_int}.

(iii)\ Noting that
$\PP(^{(r)}T\in\rmd t)=g_r(t)\rmd t$, \eqref{con_density} follows from \eqref{prop3.4} by a change of variable.

{(iv)} It can be read from \eqref{MCr} that 
\ben
\PP \big(G_{i+1} \in \rmd w \, \big|{}^{(r)}T_i = t \, \big) =
(r+i)\Theta(tw)\frac{g_{r+i+1}(t\bar w)}{g_{r+i}(t)}\rmd w\,.
\een
Comparing this with \eqref{con_density}, we obtain \eqref{gr_MC1}.
\halmos

Comparing \eqref{con_density} and \eqref{gr_MC1}, we see that, conditional on $\{{}^{(r)}T_i = t\}$, 
 $G_{i+1}$  has the same density as the first size-biased pick from $\PK^{(r+i)}(\rho \,|\, t)$.
This means that the $(i+1)$st size-biased pick from the remaining point process $\BB^{(r)}_i := \BB^{(r)} - \sum_{j=1}^i \delta_{\wt J_j}$ has the same distribution as the first size-biased pick from an independent point process $\BB^{(r+i)}$ after conditioning on their sums. 
This gives a characterisation of the sequence obtained by removing the first $k$ size-biased jumps and then renormalising it as 
\ben
\wt \bfW_k : = \bigg(   \frac{\wt J_{k+1}}{ {}^{(r)}T_k}, \frac{\wt J_{k+2}}{ {}^{(r)}T_k}, \frac{\wt J_{k+3}}{ {}^{(r)}T_k}, \cdots\bigg),
\een
as stated in the next corollary.

\begin{corollary}\label{del2}
For each $t > 0$, $r > 0$ and $k \in \N$, we have, for $w_i\in (0,1)$, $i\in\N$, 
\be\label{del1}
\PP\bigg( \frac{\wt J_{k+i}}{ {}^{(r)}T_k} \in \rmd w_i,\, i\ge 1  \, \Big| \, {}^{(r)}T_k = t    \bigg) 
=\, \PP\bigg(   \frac{\wt J_{i}}{ {}^{(r+k)}T_0}\in \rmd w_i, \, i \ge 1  \, \Big| \, {}^{(r+k)}T_0 = t    \bigg).
\ee
\end{corollary}

\medskip\noindent {\bf Proof of Corollary \ref{del2}:}\
We  compute the finite dimensional distribution for two terms. The general case is similar. Fix $t>0$, $r >0$ and $k \in \N$ and
recall \eqref{gr_MC1}. Then the lefhand side of \eqref{del1} is
\bean  
&& \PP\Big( \frac{\wt J_{k+1}}{ {}^{(r)}T_k} \in \rmd w_1,\, \frac{\wt J_{k+2}}{ {}^{(r)}T_k} \in \rmd w_2 \, \Big| \, {}^{(r)}T_k = t  \Big) \nonumber \\
 &=&
  \PP\Big( \frac{\wt J_{k+2}}{ {}^{(r)}T_k} \in \rmd w_2 \,\Big|\,  {}^{(r)}T_k = t,  \wt J_{k+1} =t w_1 \Big) \PP\Big(  \frac{\wt J_{k+1}}{ {}^{(r)}T_k} \in \rmd w_1 \, \Big| \, {}^{(r)}T_k = t  \Big) \nonumber \\
&=&
\PP\Big( \frac{\wt J_{k+2}}{{}^{(r)}T_{k+1}}  \in  \frac{ \rmd w_2}{ \bar w_1} \, \Big| \,  {}^{(r)}T_{k+1} = t \bar w_1, {}^{(r)}T_k = t \Big) \wt f_{r+k}(w_1 \, | \, t)\, \rmd w_1
\nonumber \\
&=&
 \wt f_{r+k+1} \Big( \frac{w_2}{\bar w_1} \, \Big| \, t \bar w_1\Big) \wt f_{r+k}(w_1 \, | \, t) \, \rmd w_1\, \rmd w_2\,.
\eean
Here $\bar w_1:=1-w_1$.  
We can also compute the corresponding finite dimensional distribution from the righthand side of \eqref{del1} as 
\bean
&&
 \PP\Big( \frac{\wt J_{1}}{ {}^{(r+k)}T} \in \rmd w_1,\, \frac{\wt J_{2}}{ {}^{(r+k)}T} \in \rmd w_2 \, \Big| \, {}^{(r+k)}T = t  \Big) \nonumber \\
&=&
\PP\Big( \frac{\wt J_{2}}{ {}^{(r+k)}T} \in \rmd w_2 \,\Big|\,  {}^{(r+k)}T = t,  \wt J_1 = t w_1 \Big) \PP\Big(  \frac{\wt J_{1}}{ {}^{(r+k)}T} \in \rmd w_1 \, \Big| \, {}^{(r+k)}T = t  \Big) \nonumber \\
&=&
 \PP\Big( \frac{\wt J_{2}}{{}^{(r+k)}T_1}  \in  \frac{ \rmd w_2}{ \bar w_1} \, \Big| \,  {}^{(r+k)}T_1 = t \bar w_1, {}^{(r+k)}T= t  \Big) \wt f_{r+k}(w_1 \, | \, t)\, \rmd w_1\nonumber \\
&=&
  \wt f_{r+k+1} \Big( \frac{w_2}{\bar w_1} \, \Big| \, t \bar w_1\Big) \wt f_{r+k}(w_1 \, | \, t) \, \rmd w_1\, \rmd w_2.
\eean
Comparing these proves \eqref{del1} for two terms, and analogously for $n$ terms, hence for an infinite number of terms, which may be  rearranged into decreasing order. 
This gives \eqref{del1} as stated.
\halmos

\section{Stick-Breaking Representations}\label{sub:stick}
Let $\bfV$ be defined as in \eqref{PD} with
size-biased permutation $\wt \bfV = (\wt V_1, \wt V_2, \ldots)$. Remarkable stick-breaking properties exist for $\wt \bfV$ when $\bfV$ is distributed as $\PK(\rho_\theta)$ or $\PK(\rho_\alpha)$.
Recall that for each $x > 0$,
\ben\label{levyden}
\rho_\theta(x) = \theta e^{-x}/x, \ \theta > 0, \ \text{and} \
\rho_\alpha(x) = C \alpha x^{-\alpha-1}, \  0<\alpha<1, \ C > 0.
\een

By \citet[Thm.1.2]{PPY1992},  the $n$th term in $\wt \bfV$ can be written in the product form
\ben\label{stick-break2}
\wt V_n = (1-U_n) \prod_{i=1}^{n-1} U_i, \quad n\in\N
\een
(with $\prod_{i=1}^{0}\equiv 1$).
When $\bfV$ follows a $\PK(\rho_\theta)$ distribution, the $(U_i)$ are i.i.d. Beta($\theta, 1)$ variables. When $\bfV$ is distributed as $\PK(\rho_\alpha)$, the $(U_i)$ are independent Beta($i\alpha, 1-\alpha)$ variables. 

To derive the corresponding stick-breaking representations for vectors having distributions in $\PK^{(r)}(\rho_\theta)$ and $\PK^{(r)}(\rho_\alpha)$,  recall the sequence of sums remaining after successive size-biased picks, $({}^{(r)}T_i)$, $i\in\N$, defined in \eqref{sb_T}, and denote the successive {\it residual fractions} by 
 \ben
 {}^{(r)}U_i := \frac{{}^{(r)}T_i}{{}^{(r)}{T_{i-1}}}, \ i\in\N.
 \een
Then the $n${th} term in $\wt \bfW$ can be represented as
\be\label{stick-break}
\wt W_n = (1- {}^{(r)}U_n) \prod_{i=1}^{n-1} {}^{(r)}U_i, \ n\in\N.
\ee

\begin{theorem}[Stick-Breaking]\label{stick}
Fix $r > 0$. 

{\rm (i)}\  Suppose $\bfW$ is distributed as $\PK^{(r)}(\rho_\theta)$. Then the stick-breaking factors in \eqref{stick-break} are distributed as
\be\label{ss1}
\big({}^{(r)}U_i\big)_{i\in\N} \eqd 
\big(B_i( \Gamma_r\theta, 1 )\,\big)_{i\in\N},
\ee
where, for each $v>0$, $\big(B_i(v\theta, 1))_{i\in\N}$ are i.i.d. Beta($v \theta, 1)$ variables, 
and $\Gamma_r$ is a Gamma($r,1)$ random variable, independent of the $(B_i)$.

{\rm (ii)}\ 
Suppose $\bfW$ is distributed as $\PK^{(r)}(\rho_\alpha)$. Then the corresponding $({}^{(r)}U_i$) in \eqref{stick-break} are distributed as Beta($i \alpha, 1-\alpha)$, $i\ge 1$, not depending on  $r$.
\end{theorem}

\medskip\noindent {\bf Proof of Theorem \ref{stick}:}\
We make use of the fact that $\PK^{(r)}(\rho)$ has the same law as $\PK(\Gamma_r \rho)$. Recall from \citet[Theorem 2.1]{PPY1992} that if $\bfV$ has law $\PK(v \rho_\theta)$, then the corresponding residual fractions $(U_i)$ are independent Beta$(v\theta, 1)$ variables. 
Randomising $v$ by an independent Gamma$(r,1)$ random variable, we obtain \eqref{ss1}.

Next suppose $\bfV$ has law $\PK(v \rho_\alpha)$. Then the corresponding residual fractions $(U_i)$ are independent Beta$(i\alpha, 1-\alpha)$ variables, whose distribution does not depend on $v$. Thus, randomising $v$ does not change the distribution. 
\halmos

By Theorem \ref{stick},  $\PK^{(r)}(\rho_\alpha)$ comprises the same laws as $\PK(\rho_\alpha)$ for all values of $r > 0$.
Thus, as a characteristic of $\PK(\rho_\alpha)$, 
\[ {}^{(r)}T_i \quad \text{is independent of} \quad ({}^{(r)}U_1, \ldots, {}^{(r)}U_i)\,, \quad i \in \N.
\]
However the sequence $({}^{(r)}T, {}^{(r)}T_1, \ldots )$ of remaining sums has different dynamics, as will be elucidated next. Let $\rho = \rho_\alpha$ for the rest of this section. Recall the definition of $g_{r}$ in \eqref{gd}.

\begin{proposition}\label{abs}
Fix $r > 0$ and $\rho = \rho_\alpha$.
Then 
\be\label{marginalT}
\PP\big({}^{(r)}T_n \in \rmd t \big) = L_n t^{-n\alpha} g_{r+n}(t)\rmd t, \ {\rm where}\  L_n = r^{[n]} (C\Gamma(1-\alpha))^n \frac{\Gamma(n\alpha + 1)}{\Gamma(n+1)}\,,
\ee
and thus
\be\label{abs2}
\frac{\PP\big({}^{(r)}T_{n} \in \rmd t \big)}{\PP \big({}^{(r+n)}T \in \rmd t \big)} = 
L_nt^{-n\alpha}.
\ee
\end{proposition}

\medskip\noindent {\bf Proof of Proposition \ref{abs}:}\
Fix $r > 0$. For each $n \in \N$, we first derive the joint density of $\left({}^{(r)}T_n, \, {}^{(r)}U_1, \, {}^{(r)}U_2, \ldots, {}^{(r)}U_n\right)$ by change of variables using the joint density of $\left({}^{(r)}T_0, \,{}^{(r)}T_1, \, {}^{(r)}T_2, \ldots, {}^{(r)}T_n\right)$ in \eqref{jointT0} with $\rho = \rho_\alpha$.
For simplicity, we only consider the case $n=2$. An analogous derivation holds for $n >2$.
For $t_2>0$, $0<u_i<1$, $1\le i\le 2$, the joint density of $({}^{(r)}T_2, {}^{(r)}U_1, {}^{(r)}U_2 )$ is

\ben 
h(t_2, u_1,  u_2) = f\Big(\frac{t_2}{u_1 u_2 },  \frac{t_2}{u_2},t_2 \Big) t_2^2u_1^{-2} u_2^{-3};
\een
where $f$ is defined in \eqref{jointT0} and $t_2^2u_1^{-2} u_2^{-3}$ is the 
Jacobian 
from the change of variables.  Expanding the expression in \eqref{jointT0} with $\Theta(x) = C \alpha x^{-\alpha}$, we get $h(t_2, u_1, u_2)$ equal to
\begin{align}\label{jointU4}
 &r^{[2]} g_{r+2}(t_2) u_1^{-1} u_2^{-1} \Theta\bigg(\frac{t_2}{u_1u_2} \bar u_1\bigg) \Theta\bigg(\frac{t_2}{u_2} \bar u_2\bigg) \nonumber \\
=\, & r^{[2]} g_{r+2}(t_2) (C\alpha)^2 t_2^{-2\alpha}  \big(u_2^{2\alpha-1}\bar u_2^{-\alpha}\big) \big( u_1^{\alpha-1} \bar u_1^{-\alpha} \big)
\nonumber \\
=\, & \frac{r^{[2]}}{K_{2}} \, g_{r+2}(t_2)  t_2^{-2\alpha} \cdot  \bigg[ \frac{\Gamma(1+\alpha)}{\Gamma(2\alpha) \Gamma(1-\alpha)}u_2^{2\alpha-1}\bar u_2^{-\alpha}\bigg] 
\bigg[ \frac{\Gamma(1)}{\Gamma(\alpha) \Gamma(1-\alpha)}  u_1^{\alpha-1} \bar u_1^{-\alpha} \bigg]
\end{align}
where 
\[ 
K_{2} =\frac { \prod_{i=0}^1 \Gamma(1+i\alpha)}
{(C\alpha)^{2}\Gamma^2(1-\alpha)\prod_{i=1}^2 \Gamma(i\alpha)}.
\]
Integrate \eqref{jointU4} with respect to $u_1, u_2$ to get \eqref{marginalT}.
Equality \eqref{abs2} is immediate from \eqref{marginalT}.
\halmos

The next corollary gives a moment formula for ${}^{(r+n)}T$, $n \in \N$.
				
\begin{corollary}\label{mom}
Fix $r> 0$ and keep $\rho = \rho_\alpha$. For each $n \in \N$,
\be\label{moment}
L_n^{-1} = \EE\big( {}^{(r+n)}T^{-n\alpha} \big)= \frac{\Gamma(n+1)}{\Gamma(n\alpha +1)} \frac{1}{(C\Gamma(1-\alpha))^n r^{[n]}}  =  \frac{1}{r^{[n]}}\EE\big( T^{-n\alpha} \big)  \,. 
\ee
\end{corollary}

\medskip\noindent {\bf Proof of Corollary \ref{mom}:}\  
For each $n \in \N$,  integrate \eqref{marginalT} to get 
\[
1 = \int_0^\infty L_n \, g_{r+n}(t_n) \,  t_n^{-n\alpha}\, \rmd t_n = L_n \, \EE\big({}^{(r+n)}T^{-n\alpha}\big),
\]
with $L_n = r^{[n]}\, (C\Gamma(1-\alpha))^n\, \Gamma(n\alpha + 1)/\Gamma(n + 1)$. This gives the first and second equalities in \eqref{moment}. 
From \citet[Eq. (2.n)]{PPY1992}, 
\[\EE(T^{-n\alpha}) = \frac{\Gamma(n+1)}{\Gamma(n\alpha+1) (C\Gamma(1-\alpha))^n}.
\]
Comparing this with the value for $L_n^{-1}$, we get the last equality in \eqref{moment}.
\halmos

\section{$\PD_\alpha^{(r)}$ Arising from a Trimmed Stable Process}\label{sub:trimmed}

Recall the distribution $\PD_\alpha^{(r)}$ derived from the trimmed $\alpha$-stable subordinator in \eqref{PDr}.
Let $\Delta_1 > \Delta_2 > \ldots$ be the ordered jumps of an $\alpha$-stable subordinator $(S_t, 0<t<1)$. Then $\sum_i \delta_{\Delta_i}$ forms a Poisson process with intensity measure $\Lambda(\rmd x) := \rho_\alpha(x)\rmd x$.
For $r \in \N$, denoting the $r$-trimmed process up till time 1 by ${}^{(r)}S_1 = S_1 - \sum_{i=1}^r \Delta_i$, we have
\be\label{PD_r}
 ( V_n^{(r)}, n\ge 1) : =
\bigg(\frac{\Delta_{r+1}}{{}^{(r)}S_1}, \frac{\Delta_{r+2}}{{}^{(r)}S_1}, \frac{\Delta_{r+3}}{{}^{(r)}S_1}, \ldots \bigg)
 \sim \PD_\alpha^{(r)}.
\ee

In this section we derive some analogous properties for $\PD_\alpha^{(r)}$.
First note that
$\PD_\alpha^{(r)}$ has the same law as $\PK^{(r)}(\rho_\alpha^*)$, for $\rho_\alpha^* (x):= \alpha x^{-\alpha-1}{\bf 1}_{0<x < 1}$.
To see this, we know from \citet[Lemma 24]{PY1997}, that, conditionally on $\Delta_{r}$, the point process 
\be\label{nj}
\BB^T := \sum_{i\in\N} \delta_{J_{(i)}},  \ \text{where } J_{(i)} = \Delta_{r+i}/\Delta_{r},
\ee
is a Poisson process with intensity measure $(\Delta_{r})^{-\alpha} \Lambda(\rmd x){\bf 1}_{\{0<x<1\}}$. Since $(\Delta_{r})^{-\alpha} \eqd \Gamma_r/C$, $\BB^T$ is a negative binomial point process, $\BN(r, \rho_\alpha^*)$.  
Thus the vector 
\[
\bigg(\frac{J_{(1)}}{{}^{(r)}S_1/ \Delta_{r}}, \frac{J_{(2)}}{{}^{(r)}S_1/ \Delta_{r}}, \frac{J_{(3)}}{{}^{(r)}S_1/\Delta_{r}}, \ldots\bigg),
\]
which is equal to the  lefthand side  of \eqref{PD_r}, has distribution $\PK^{(r)}(\rho_\alpha^*)$.

Due to the restriction to $(0,1)$ in $\rho_\alpha^*$,   $\PD_\alpha^{(r)}$ is distinct from $\PD(\alpha, 0)$
(unlike for $\PK^{(r)}(\rho_\alpha)$). The next theorem gives a characterisation of the stick breaking sequences as in \eqref{stick-break}.
\begin{theorem}\label{stick2}
The joint distribution  of ${}^{(r)} T_n$ and
${}^{(r)}U_1, \, {}^{(r)}U_2, \ldots,  {}^{(r)}U_n$   for $\PD_\alpha^{(r)}$
can be written as 
\be\label{h6}
\big({}^{(r)}T_n, \, {}^{(r)}U_1, \, {}^{(r)}U_2, \ldots,  {}^{(r)}U_n\big)
\eqd \big( Y_{d(U_1, \ldots, U_n)},  U_1,  U_2, \ldots,  U_n\big)�,
\ee
where the  $(U_i)$ are independent Beta$(i\alpha, 1-\alpha)$ rvs,
\[d(u_1, \ldots, u_n):= \min_{1\le i \le n} \prod_{j=i}^n u_j/\bar u_i,
\] 
with 
 $\bar u_i := 1 - u_i$,  and for each $c > 0$, 
$Y_c := ({}^{(r+n)}T)^{-n\alpha}{\bf 1}_{\{{}^{(r+n)}T < c\}}$.
\end{theorem}

\medskip\noindent {\bf Proof of  Theorem \ref{stick2}:}\
Let $r > 0$ and $\rho = \rho_\alpha^*$.
Similar to \eqref{jointU4}, we can derive for $t_n>0$, $0<u_i<1$, $1\le i\le n$,  and $n \in \NN$, the joint density of 
$\left({}^{(r)}T_0, \,{}^{(r)}T_1, \, {}^{(r)}T_2, \ldots, {}^{(r)}T_n\right)$ as
\begin{align}\label{jointUTn}
h(t_n, u_1, \ldots, u_n) 
&=
\frac{r^{[n]}}{K_{n}}  g^*_{r+n} (t_n)  t_n^{-n\alpha} \nonumber \\
&\qquad \times 
 \prod_{i=1}^n \frac{\Gamma(i\alpha+1-\alpha)}{\Gamma(i\alpha)\Gamma(1-\alpha)} u_i^{i\alpha-1} \bar u_i^{-\alpha} {\bf 1}_{\{ t_n <\prod_{j=i}^n u_j/\bar u_i\}}\nonumber \\
 &= \frac{r^{[n]}}{K_{n}} t_n^{-n\alpha} g^*_{r+n}(t_n) {\bf 1}_{\{ t_n <d(u_1, \ldots, u_n)\}} \times
 \prod_{i=1}^n \beta_{i\alpha, 1-\alpha}(u_i).
\end{align}
Here $g^*_r$ is the corresponding density of the sum of points in $\BN(r, \rho_\alpha^*)$, 
$\beta_{a,b}$ is the density of a Beta(a,b) distribution, and
\ben\label{Kn}
 K_{n}
= 
\frac{\prod_{i=0}^{n-1} \Gamma(1+i\alpha)}
{\alpha^{n} \Gamma^n(1-\alpha) \prod_{i=1}^n \Gamma(i\alpha)} 
=
\frac {\Gamma(n+1)}{\Gamma^n(1-\alpha) \Gamma(n\alpha+1)}.
\een
The indicator function  in \eqref{jointUTn} reflects the restriction of $x$ to the interval $(0,1)$ in $\rho_\alpha^*$.
\halmos

\begin{remark}\label{rem3d} 
{\rm 
(i)\, By integrating \eqref{jointUTn}, we get the identity 
\begin{align*}\label{depen}
K_{n} &= r^{[n]}\int_{u_1=0}^1 \cdots \int_{u_n=0}^1 \int_{t_n = 0}^{d(u_1, \ldots, u_n)} 
t_n^{-n\alpha} g^*_{r+n}(t_n) \rmd t_n \nonumber \\
&\hspace{2in} \times \prod_{i=1}^n \beta_{i\alpha, 1-\alpha}(u_i)\rmd u_1 \cdots \rmd u_n.
\end{align*}

(ii)\ 
For a stick-breaking representation, as in \eqref{stick-break}, the size-biased permutation of $(V_n^{(r)})$, denoted by $(\wt V_n^{(r)})$,  can be written as
\be\label{wVsb}
\wt V_n^{(r)} = (1- {}^{(r)}U_n) \prod_{i=1}^{n-1} {}^{(r)}U_i.
\ee
The joint distribution of $({}^{(r)}U_i)_{1\le i\le n}$ can be computed from \eqref{h6}, in which we note that $U_1,  U_2, \ldots,  U_n$ are individually independent, 
but dependence overall is introduced via the connection with the $Y$ term.
In this respect the result is different from the $\PD(\alpha,0)$  situation, as we would expect, but the distribution of $\wt V_n^{(r)}$ as given by \eqref{wVsb} is sufficiently  explicit to enable computations or simulations.

(iii)\ Although motivated by the idea of trimming an integer number $r$ of large jumps, our formulae once derived are valid for $r>0$, and available for modelling purposes in this generality.

(iv)\ We may set $r=0$ in \eqref{PD_r} to have the distribution of $\PD_\alpha^{(r)}$ reduce to that of $\PD(\alpha,0)$. But we cannot take $r=0$ in \eqref{nj} 
with the idea that the size-biased distribution associated with 
$\PD_\alpha^{(r)}$ might then reduce to the one associated with $\PD(\alpha,0)$.
Note that $\BB^T$ is not defined for $r=0$ (its points $\Delta_{r+i}/\Delta_i$ are not defined for $r=0$).  
Setting $r=0$ in \eqref{jointUTn}, which results from an analysis of $\BB^T$,
is not permissible.
}
\end{remark}

By restricting $r$ to be an integer, we can further construct a vector $(V_n^{(r)})$ from independent beta random variables and characterise the law of the sequences in $\PD_\alpha^{(r)}$ as a shifted version of $\PD(\alpha, 0)$ which can in turn
be characterised by a change of measure formula.

\begin{lemma}\label{changeM}
{\rm (i)}\, Let $(R_i, i\ge 1)$ be a sequence of independent Beta$((r+i) \alpha, 1)$ variables. Define 
\[
V_n^{(r)} := \frac{\prod_{i=1}^{n-1}R_i}{1+R_1+R_1R_2 + \cdots} =  Y_n \prod_{i=1}^{n-1}(1-Y_i),
\] with stick-breaking factors $Y_n = (1+ \Sigma_n)^{-1}$ and $\Sigma_n = R_n+R_nR_{n+1}+ \cdots$. Then $(V_n^{(r)}, n\ge 1)$ has law $\PD_\alpha^{(r)}$.

{\rm (ii)}\, Let $(V_i^{(r)})$ be distributed as $\PD_\alpha^{(r)}$ and $(V_i)$ be distributed as $\PD(\alpha, 0)$. Then for any nonnegative measurable function $f$, we have 
\[
\EE \{ f(V_1^{(r)}, V_2^{(r)}, \ldots) \} = \EE \Big\{ \frac{\EEEE_1^r}{r!}f(V_1 , V_2, \ldots)\Big\},
\]
where $\EEEE_1 = \lim_{n\to \infty} nV_n^\alpha/V_1^\alpha$, and the limit holds almost surely and in $p$th mean for all $p \ge 1$.
\end{lemma}

\noindent{\bf Proof of Lemma \ref{changeM}:}
(i) Recall that  $(\Delta_1, \Delta_2, \ldots)$ comprise the points of a Poisson process with intensity measure $\Lambda(\rmd x) = \rho_\alpha(x) \rmd x$. 
Write $R_n = \Delta_{r+n+1}/\Delta_{r+n}$. Then by \citet[Prop.8]{PY1997}, the sequence of successive ratios $(R_1, R_2, \ldots )$ is of  independent Beta$({(r+n)\alpha, 1})$ variables.

Since the sequence 
\[
\Big(\frac{\Delta_{r+1}}{\sum_{i\ge r+1}\Delta_i}, \, \frac{\Delta_{r+2}}{\sum_{i\ge r+1}\Delta_i}, \ldots \Big) \sim \PD_\alpha^{(r)},
\]
the $n$th term can be expressed as 
\begin{align*}
\frac{\Delta_{r+n}/\Delta_{r+1}}
{
\Delta_{r+1}/\Delta_{r+1}+ \Delta_{r+2}/\Delta_{r+1}+\cdots
} 
&=  \frac{\prod_{i=1}^{n-1}R_i}{1 + R_1 + R_1R_2+\ldots} \nonumber \\
&=  \frac{\Delta_{r+n}}{\sum_{i\ge n}\Delta_{r+i}}\cdot  \frac{\sum_{i\ge n}\Delta_{r+i}}{\sum_{i\ge n-1}\Delta_{r+i}} \cdots  \frac{\sum_{i\ge 2}\Delta_{r+i}}{\sum_{i\ge 1}\Delta_{r+i}} \nonumber \\
&= Y_n \prod_{i=1}^{n-1} (1-Y_i)
\end{align*}
where $Y_n = \Delta_{r+n}/\sum_{i\ge n} \Delta_{r+i} = (1+\Sigma_n)^{-1}$.

(ii)
First consider the homogeneous Poisson point process $\sum_{i}\delta_{\Gamma_i}$, where $\Gamma_n := \sum_{i=1}^n E_i$ with $(E_i, i\ge 1)$ independent unit exponential random variables.
Then for any nonnegative measurable function $f$, 
the shifted sequence $(\Gamma_{r+i}, i\ge 1)$ can be characterised by a change of measure from the original process as 
\[
\EE \{f (\Gamma_{r+1}, \Gamma_{r+2}, \ldots )\} = \EE \Big\{ \frac{\Gamma_1^r}{r!} f (\Gamma_{1}, \Gamma_{2}, \ldots )\Big\}.
\]

Since for each $i\ge 1$, $\Delta_i \eqd \laminv(\Gamma_i)$, where $\laminv(x) = C^{1/\alpha}x^{-1/\alpha}$,
then for each nonnegative measurable function $f$, we again have 
\[
\EE \{f (\Delta_{r+1}, \Delta_{r+2}, \ldots )\} = \EE \Big\{ \frac{\EEEE_1^r}{r!} f (\Delta_{1}, \Delta_{2}, \ldots )\Big\}
\]
where $\EEEE_1  = \Lambda(\Delta_1, \infty)$.
Thus, normalising the jumps,  we still get
\[
\EE \{f (V_1^{(r)}, V_2^{(r)}, \ldots )\} = \EE \Big\{ \frac{\EEEE_1^r}{r!} f (V_1,V_2, \ldots ) \Big\}
\]
where $\EEEE_1 := \Lambda(\Delta_1, \infty) = C \Delta_1^{-\alpha}$. 
To write $C\Delta_1^{\alpha}$ as a function of $(V_1, V_2, \ldots )$, we note that 
by \citet[Prop. 10]{PY1997}, $\lim_{n\to \infty} n V_n^{\alpha} = C S_1$ almost surely and in $p$th mean for $p \ge 1$. Thus, 
we can write 
\[
\EEEE_1 = C \Delta_1^{-\alpha} = \frac{C}{V_1^\alpha S_1^\alpha} = \frac{\lim_{n \to \infty} n V_n^\alpha}{V_1^\alpha},
\]
concluding the proof of Lemma \ref{changeM}.
\halmos

\section{Discussion: Applications}\label{dis}

We mention some possible applications of our results.
A common tool in linguistics studies is the ``Zipf plot":  a plot of log frequencies of words, against their log-ranks.
 \cite{ggj11} 
show such  plots for some word counts in the Penn Wall St. journal. In their Figure 4, half a dozen or so of the most frequent words appear as outliers, while the rest conform closely to a $\PD(\alpha, 0 )$ fit. This suggests that  a $\PD_\alpha^{(r)}$ distribution might provide a better fit.

A similar situation occurs in  \cite{sos15}, who shows ``capital distribution curves" 
(a log plot of normalized stock capitalizations ranked in descending order, against their log-ranks) for over 20 countries listed on the NASDAQ stock exchange. 
The curves  appear to be very well fitted by a $\PD(\alpha,0)$ distribution over much of their range, but with a small number of the largest stocks as outliers -- as we might expect from this kind of data.

Known difficulties arise in fitting the general 2-parameter $\PD(\alpha,\theta)$  distribution to data; the maximum likelihood estimator of $\theta$ is inconsistent (\cite[Lemma~5.7]{Carlton1999}). Introducing the extra parameter $r$ in $\PK^{(r)}(\rho)$ may help to improve estimation of $\theta$, as well as allowing extra flexibility in data description.

In general, we expect that our generalised $\PK^{(r)}(\cdot)$ distribution could be used to extend analyses which are implicitly based on underlying Poisson point processes, to negative binomial point processes, and thereby reveal interesting features of data. More research along these lines would certainly prove profitable.

\medskip\noindent{\bf Acknowledgement}\ 
We appreciate  the referee's close reading of the paper, and helpful comments.

\renewcommand{\bibfont}{\small}
\bibliography{Library_Levy_Dec2017}
\bibliographystyle{newapa}

\end{document}